\documentclass[a4paper,12pt]{article}
\usepackage[T2A]{fontenc}
\usepackage[utf8]{inputenc}
\usepackage{amsmath,amssymb}
\usepackage[pdftex]{graphicx}
\graphicspath{{noiseimages/}}

\date{}

\vspace{5mm}

\title{\bf Weak ergodic theorem for Markov chains without invariant countably additive measures}

\author{\bf {Alexander Zhdanok} \footnote{Institute for Information Transmission Problems of the Russian Academy of Science, Moscow;} 
\footnote{e-mail: zhdanok@inbox.ru}
 \footnote{{\bf Acknowledgments:} This research was funded by the Russian Foundation of Basic Research, RFBR project number 20-01-00575-a.}} 

\begin{document}

\maketitle

\begin{abstract}
	
In this paper, we study Markov chains (MC) on topological spaces within the framework of the operator approach. We extend the Markov operator from the space of countably additive measures to the space of finitely additive measures. Cesaro means for a Markov sequence of measures and their asymptotic behavior in the weak topology are considered.
It is proved ergodic theorem that in order for the Cesaro means to converge weakly to some bounded regular finitely additive (or countably additive) measure it is necessary and sufficient that all invariant finitely additive measures are not separable from the limit measure in the weak topology. Moreover, the limit measure may not be invariant for a MC, and may not be countably additive. The corresponding example is given and studied in detail. \\



\end{abstract}



\section*{\bf Introduction}



In this paper, time-homogeneous Markov chains (MCs) on a normal topological space are studied in the framework of the operator-theoretical approach. Such an approach of the general MC study was proposed by N. Kryloff and N. Bogoliouboff (1937), which was explicitly developed then by K. Yoshida and S. Kakutani in~\cite{YoKa:1}.
MCs are given by a transition functions $p(x, E), x \in X, E \in \Sigma$, where $ X $ is some set (space) and $\Sigma$ is some sigma-algebra of subsets in $X$. The transition functions $p(x,E)$ are assumed as a countably additive measures by the second argument, i.e., classical Markov chains are considered.


The transition function $p(x,E)$ of a Markov chain defines two dual integral Markov linear bounded positive operators $ T $ and
$ A $ in spaces of measurable functions and in spaces of measures, respectively.
The Markov operators of a Markov chain extend from the space of countably additive measures to the spaces of finitely additive measures.
Such operators become topologically conjugate to Markov operators defined on spaces of measurable functions (here all measures and functions are assumed to be bounded). This construction of operators allows a wider use of the methods of functional analysis. The spaces of bounded continuous functions and the weak topologies generated by them in the spaces of measures are also involved. The fact is used that such weak topologies have inseparable sets in the space of finitely additive measures, and there appear sets of "sticky points (measures)" in this topology. The introduction of regular measures makes it possible to describe in more detail some properties of such sets. In particular, we use the notion of regularization of finitely additive irregular measures introduced and studied by us earlier in ~\cite{Zhd84} and~\cite{Zhd2003a}.

In sections 1 and 2,  we give some necessary for us concepts, facts and symbols that are rarely used in classical Probability Theory.


In section 3, an ergodic theorem is formulated for an arbitrary MC on an arbitrary normal topological space. This new theorem proves that for the weak convergence of Cesàro means from a Markov sequence of measures to some regular finitely additive measure (condition $(i)$) it is necessary and sufficient that all invariant finitely additive (irregular) measures of such an MC in the weak topology were inseparable from the above limit measure (condition $(ii)$).
In the formulation of this theorem and in its proof, it is not required that the weak limit measure for Cesàro means in the condition $(i)$ be countably additive or invariant (although this is allowed as a special case). The Markov chain is also not assumed to be Fellerian. 

In this article, we use the methods and results of our papers ~\cite{Zhd2003a}, ~\cite{Zhd2003b} and ~\cite{Zhd84}. The proved weak ergodic theorem is a generalization of Theorem 4 from ~\cite{Zhd1981} and Theorem 13.1 from ~\cite{Zhd2003b}. 


At the end of the article, in Section 4, we construct an example of an MC, on which we demonstrate in detail the application of Theorem 3.1. 


\section{Finitely additive measures} 

Let $X$ be an arbitrary set and let $\Sigma$ be an $\sigma$-algebra (or algebra) of its subsets.

We assume $\Sigma$ to contain all singletons of $X$.

If $X$ is a topological space with topology $\tau=\tau_{X}$ then $\mathfrak{A}=\mathfrak{A}_{X}=\mathfrak{A}_{\tau}$ are the Borel algebra and  $\mathfrak{B}=\mathfrak{B}_{X}=\mathfrak{B}_{\tau}$ are the Borel $\sigma$-algebra on $X$ generated by $\tau$. Assume that $X$ is normal (all metric space are normal). \medskip
           
Following Dunford and Schwartz notations~\cite{DS1958}, we denote by 

$ba(X,\Sigma)$ the Banach space of all bounded finitely additive measures $\mu: \Sigma \to R $ with  norm the total variation of a measure on $X$ ($\|\mu\|=Var(\mu , X)$), and by 

$ca(X,\Sigma)$ the Banach space of all bounded countably additive measures $\mu: \Sigma \to R \; $ also with total variation as the norm. 

Finitely additive measures are also referred to as charges in the literature.\medskip

{\bf Definition 1.1.} (Yosida and Hewitt,~\cite{YoHew1952}). A nonnegative finitely additive measure $\mu \in ba(X, \Sigma)$ is called  purely finitely additive if for every countable additive measure $\lambda \in ca(X, \Sigma)$, $0 \le \lambda \le \mu$ implies $\lambda = 0$. A measure $\mu \in ba(X, \Sigma)$ is called purely finitely additive if both nonnegative measures $\mu^{+}$ and $\mu^{-}$ of its Jordan decomposition $\mu = \mu^{+}-\mu^{-}$ are purely finitely additive.

{\bf Theorem 1.1.} (Yosida and Hewitt,~\cite{YoHew1952}) {\it Each finitely additive measure $\mu \in ba(X, \Sigma)$ is uniquely representable as $\mu=\mu_{1}+\mu_{2}$, where $\mu_{1} \in ca(X, \Sigma)$ is a countably additive measure and $\mu_{2} \in ba(X, \Sigma)$ is a purely finitely additive measure.}

Like countably additive measures, purely finitely additive measures form a vector Banach subspace in $ba(X, \Sigma)$ which we denote by $ pfa(X, \Sigma)$. Theorem 1.1 can be treated as an assertion on direct decomposition of the measure space:
                    $$ ba(X, \Sigma)=ca(X, \Sigma) \oplus  pfa(X, \Sigma).$$

A purely finitely additive measure vanishes on every finite set. \medskip   

Given an arbitrary $(X, \Sigma)$ with $\sigma$-algebra $\Sigma$, denote by ${ B(X, \Sigma)}$ the Banach space of all bounded ${\bf \Sigma}$-measurable functions $f: X \to R$ with the $\sup$-norm.\medskip

{\bf Definition 1.2}~(\cite{DS1958}, Chapter III). Let $X$ be a normal topological space and $\Sigma$ an arbitrary algebra in  $X$. A set $E \in \Sigma$ is {\it regular} for a measure $\lambda \in ba(X, \Sigma)$ if, for	every $\varepsilon >0$, there exist $F,G \in \Sigma$ such that $\bar{F} \subset E \subset \displaystyle{\buildrel\circ\over G}$ and $Var(\mu, G \setminus F)< \varepsilon$ (here $\bar{F}$ is the closure of $F$ and $\displaystyle{\buildrel\circ\over G}$ is the interior of G). A measure $\lambda \in ba(X, \Sigma)$ is called {\it regular} if each $E \in \Sigma$  is regular.\medskip

In this definition, the algebra $\Sigma$ is not necessarily Borel in $X$, i.e., $\Sigma \neq \mathfrak{A}$ is possible. This is used in further constructions in~\cite{DS1958}.

We use the standard notations from~\cite{DS1958}: 

${rba(X, \Sigma)}$ is the Banach space of all regular finitely additive bounded measures 
 on  $(X, \Sigma)$; 
 
${ rca(X, \Sigma)}$ is the Banach space of all regular countably additive bounded measures on $(X, \Sigma)$. 

It is known that, in a metric space $X$, if $\Sigma=\mathfrak{B}$ then every countably additive measure is regular, i.e., $ca(X, \mathfrak{B}) = rca(X, \mathfrak{B})$. \medskip

Let $X$ be a normal topological space. We denote by ${C(X)}$ the Banach space of all bounded continuous functions $ f: X \to R$ with $sup$-norm. It is known that $C(X) \subset B(X, \mathfrak B)$.

Recall that there is a topological duality between the vector Banach spaces of functions and measures (see~\cite{DS1958}, Chapter IV). The conjugate spaces are:
$$B^{*}(X,\Sigma)=ba(X,\Sigma) \text{ for an arbitrary } (X,\Sigma),$$ 
$$C^{*}(X)= rba(X, \mathfrak A) \text{ for a normal topological } X,$$ 
$$C^{*}(X)=rca(X, \mathfrak B) \text{ for a Hausdorff compact space } X.$$ 

The equalities denote isometric isomorphisms.

Let $M$ be an arbitrary space of measures. We use the notations $S_{M} = \{ \mu \in M: \mu \ge 0, \mu(X) = 1 \}$. Thus, $S_{ca}$ is the set of all traditional countably additive probability measures on $(X,\Sigma)$. All measures on sets $S_{M}$  for arbitrary $M$ also will be called probability measures.

If $X$ is a topological space then, in many problems, the original measure $\lambda \in ba(X, \Sigma)$ can be replaced by a regular measure $\overline{\lambda}$ "stuck" to it in the topology $\tau_{C}$ generated by $C(X)$ in $ba(X, \Sigma)$. Such a procedure was studied in detail by the author in~\cite{Zhd84}  and given also in ~\cite{Zhd2003a}. In this connection, we now recall some facts that are necessary for the exposition.

{\bf Theorem 1.2} (\cite{Zhd84}, Theorem 1, and~\cite{Zhd2003a}, Theorem 1.3). Suppose that $X$ is normal. Then, for every $\lambda \in ba(X, \Sigma)$, there exists a unique $\overline{\lambda} \in rba(X, \mathfrak{A})$ such that $\int f d\lambda = \int f d{\overline{\lambda}}$ for every $f \in C(X)$. Moreover, $\lambda(X) = \overline{\lambda} (X)$; if $ \lambda\ge 0$ then $ \overline{\lambda} \ge 0$; if $\lambda \in ca(X, \mathfrak{B})$ then $\overline{\lambda} \in rca(X, \mathfrak{B})$ (as the extension of $\overline{\lambda} \in rca(X, \mathfrak{A})$ to $\mathfrak B$).

{\bf Definition 1.3} (\cite{Zhd84} and~\cite{Zhd2003a}). Given $\lambda \in ba(X, \mathfrak A)$, we call the measure $\overline{\lambda} \in rba(X, \mathfrak{A})$ corresponding to $\lambda $ by Theorem 1.2 the {\it regularization} of $\lambda$.

{\bf Corollary 1.1.} (\cite{Zhd84} and~\cite{Zhd2003a}). If $X$ is a Hausdorff compact space then, for every $\lambda \in ba(X, \mathfrak A)$, its regularization $\overline{\lambda}$ belongs to $rca(X, \mathfrak{B})$.

{\bf Definition 1.4} (\cite{Zhd84} and~\cite{Zhd2003a}). Assume that $\mu \in rba(X, \mathfrak A)$ and $\mu \ge 0$. The set $\mathfrak R \{ \mu \} = \{ \lambda \in ba(X, \mathfrak A) : \lambda \ge 0, \overline{\lambda} = \mu \}$ is called the class of {\it C-equivalent measures for $\mu$}.

{\bf Theorem 1.3} (\cite{Zhd84}, Theorem 5, and~\cite{Zhd2003a}, Theorem 1.4). Let $\mu \in rba(X, \mathfrak A)$ and $\mu \ge 0$. The set $\mathfrak R \{ \mu \}$ is convex and compact in the $\tau_{B}$-topology of $ba(X, \mathfrak B)$ ($\tau_{B}$ is the *-weak topology on $ba(X, \mathfrak B)$ generated by the pre-conjugate space $B(X, \mathfrak B)$).

It should be noted that, as a matter of fact, this natural pair of measures ($\lambda, \overline{\lambda}$) was used by some other authors as an intermediate technical tool (without studying the interrelation between $\lambda$ and $\overline{\lambda}$ in detail).


\section{Markov operators} 

{\bf Definition 2.1.}  A transition function (transition probability) $p(x,E)$ on a measure space $(X, \Sigma)$ is a mapping $p: X \times \Sigma \to [0,1]$ satisfying the conventional conditions
$$
p( \cdot , E)  \in B(X,\Sigma), \; \forall E \in  \Sigma;
$$
$$ p( x, \cdot ) \in  ca(X,\Sigma), \;  \forall x \in  X ;$$
$$ p( x, X) = 1, \; \forall x \in X.
$$ 

The transition function $p(x,E)$ is countably additive for the second argument, i.e. classical.\medskip

{\bf Definition 2.2.} By the {Markov operators}, we mean the two operators $T$ and $A$ defined explicitly as follows:
$$ 
T: B(X,\Sigma) \to B(X,\Sigma),\;\;\; (Tf)(x)=Tf(x) \stackrel{\mathrm{def}}{=} \int_{X}{f(y)p(x,dy)},
$$
$$
\text{where } f \in B(X,\Sigma), \; x \in X; 
$$
$$ 
A: ca(X,\Sigma) \to ca(X,\Sigma),\;\;\; (A \mu)(E)=A \mu(E) \stackrel{\mathrm{def}}{=} \int_{X}{ p(X,\Sigma) \mu(dx)},$$ 
$$ 
\text{\rm where }\;\; \mu \in ca(X,\Sigma), E \in \Sigma.
$$  


Assume that $\mu_{0} \in S_{ca} \; \text{\rm and} \; \mu_{n}=A^{n}\mu_{0}=A\mu_{n-1}, \; n=1,2, \ldots$. An MC can be identified with the sequence of probability measures $\{\mu_{n}\}=\{\mu_{n}(\mu_{0})\}$ depending on the initial measure $\mu_{0}$ as a parameter. Therefore, every MC can be regarded as an iterative process generated by a positive linear operator $A$ on a space of measures. The sequence $\{\mu_{n}\}$ will be called a Markov sequence of measures. \medskip

It is known that {\it for every countably additive MC the Markov operator $A$ of Definition 2.2 is uniquely extendable from $ca(X,\Sigma)$ to a linear operator $\tilde A$ on $ba(X,\Sigma)$, preserving positivity, isometry on the cone, boundedness, the norm, and explicit form
$$
\tilde A: ba(X,\Sigma) \to ba(X,\Sigma), \;\;$$
$$ (\tilde A \mu)(E) \stackrel{\mathrm{def}}{=} \int_{X}{p(x, E) \mu(dx),}
\;\; \mu \in ba(X,\Sigma), E \in \Sigma.
$$ 
Moreover, $\tilde A$ is topologically adjoint to the operator $T$ of Definition 2.2, i.e., $T^{*}=\tilde A$ with $B^{*}(X,\Sigma)=ba(X,\Sigma)$.} \medskip

We call the extension $\tilde A$ of the Markov operator $A$ the {finitely additive extension of $A$}. Like $A$, we call $\tilde A$ a {\it Markov operator}.
Below, we often identify $\tilde A$ and $A$ without specifying their domains of definition.

Suppose that $\mu_{0} \in ba(X,\Sigma)$ is such that $\mu_{0} \ge 0$ and $\|\mu_{0}\|=\mu_{0}(X)=1, \; \text{\rm i.e.} \; \mu_{0} \in S_{ba}$. Then $\tilde A$ generates the sequence of finitely additive measures $\mu_{n}= \tilde A \mu_{n-1}= {\tilde A}^{n} \mu_{0} \in ba(X,\Sigma), \; n=1,2,\ldots$.    Following our ideology, such an iterative process can be treated as a countably additive MC extended to the space of finitely additive measures.

Emphasize that we carry out a finitely additive extension of the operator A and MC itself for a transition probability, which is still countably additive.

At the same time, it is possible to consider MC with a finitely additive transition probability. Such MCs are called finitely additive MCs and are studied, for example, in ~\cite{Zhd2003a},~\cite{Rama1981} and~\cite{ZhdKh2022}. We do not consider such MCs in this paper.

We will need the following two \v{S}idak theorems.

{\bf Theorem 2.1} (\v{S}idak, 1962, \cite{Si1962}). {\it For every MC on an arbitrary measure space \; $(X,\Sigma)$, there exists an invariant finitely additive measure
$ 
\lambda \in ba(X,\Sigma), \lambda\ge0, \; \lambda (X)=1, \; 
\; \text{\rm i.e.}\; \lambda \in S_{ba}, \; \lambda =\tilde A \lambda \; \text{\rm and}\;
$
$$
\lambda(E)=\int{p(x,E) \lambda (dx)}, \; \forall E \in \Sigma.
$$ }

{\bf Theorem 2.2} (\v{S}idak, 1962, \cite{Si1962}). {\it Suppose that we have $\lambda=\tilde A\lambda $ for an arbitrary MC and some $\lambda \in S_{ba}$. If $\lambda=\lambda_{1}+\lambda_{2}$ is the decomposition of $\lambda$ into the sum of a countably additive and puy finitely additive measures then $\lambda_{1}=\tilde A \lambda_{1}$ and $\lambda_{2}=\tilde A \lambda_{2}$.}
\medskip

It follows from this Theorem 2.2 that in many cases it is sufficient to consider countably additive or finitely additive invariant measures separately without using their linear combinations.

For each space of measures $M$ we use in this work, denote the set of positive normalized invariant measures for the Markov operator $\tilde A$ by $\Delta_{M}=\{\mu \in S_{M}: \mu=\tilde A \mu \}$. For $\Delta_{ba}$, we sometimes omit the index: $\Delta=\Delta_{ba}$. In particular, $\Delta{pfa}=\{\mu \in \Delta: \mu \text{\rm~-- is purely finitely additive}\}$.

In Theorem 2.1 it is asserted that $\Delta_{ba} \neq \varnothing$ for any MC. However, it is possible that $\Delta_{ca} = \varnothing$ or $\Delta_{pfa} = \varnothing$.\medskip

Note 2.1. 
A systematic study of MCs with operators extended to the space of finitely additive measures was carried out in numerous works by S. R. Fogel and his colleagues (see, for example, \cite{Fog1966}, 1966, and for a more detailed bibliography see~\cite{Zhd2003a} and~\cite{Zhd2003b}). The MC phase space in these papers is topological and the MC transition probability is countably additive. In these papers, in particular, analogues of \v{S}idak theorems given above are proved for a particular case of a topological phase space. In the main constructions in Fogel papers, MCs are assumed to be Fellerian and weak topologies in measure spaces are also used.

Our papers ~\cite{Zhd2003a} and ~\cite{Zhd2003b} show points of agreement with the results of Fogel and his colleagues. Theorem 3.1, proved by us below, resembles the beginning of Theorem 1 from \cite{Fog1966}, since both Cesàro means for MC are considered and regular finitely additive measures are used. However, the statements of these two theorems are essentially different.
           
   
\section{Main rezult: Weak Ergodic Theorem}

Convergence of an MC in the $\tau_{C}$-topology, i.e., weak convergence in the probabilistic terminology, is closely connected with invariant purely finitely additive measures. We present some of our results concerning this matter. The main peculiarity of Theorem 3.1 is that {\it we do not presuppose existence of an invariant countably additive (i.e., classical "probability") measure for the MC}. 

We denote the {Cesaro means} from the Markov sequence of measures for a initial measure $\eta \in S_{ba}$ as follows:
$$\lambda_n=\lambda_n^\eta = \frac{1}{n}  \sum\limits_{k=1}^n {{A^{k}} \eta},\;\; n \in N.$$  

{\bf{Theorem} 3.1.} {\it Suppose that $X$ is a normal topological space. We have an arbitrary MC on $(X ,\mathfrak{B})$, where $\mathfrak{B}$ is the Borel sigma-algebra. Let $\mu \in rba(X,\mathfrak{A})$, $\mu \in S_{rba}$ is a some fixed regular finitely additive probability measure defined on the Borel algebra $\mathfrak{A}$. Then the following two statements (conditions) are equivalent:
	
$$ (i) \;\;\;
\int_{X}{f(x)\lambda_{n}^{\eta}(dx)} \to \int_{X}{f(x)\mu(dx)},\;\; {\text{when}} \;\; n\to \infty,$$
for any initial finitely additive measure $\eta \in ba(X,\mathfrak{B})$, $\eta \in S_{ba}$ and for any continuous function $ f \in C(X)$;

$$(ii)\;\;\; \int_{X}{f d\zeta}=\int_{X}{f d\mu},$$
for any invariant finitely additive measure $\zeta \in ba(X,\mathfrak{B})$, $\zeta \in S_{ba}$, $\zeta \in \Delta_{ba}$,
and for any continuous function $ f \in C(X)$.}

Comment 3.1. Condition $(ii)$ means that for all measures $\zeta \in \Delta_{ba}$ their regularization $\bar \zeta = \mu$, or, in other terms, $\Delta_{ba} \subset \mathfrak R \{ \mu \}$, where $\mathfrak R \{ \mu \}$ is the class $\tau_{C}$-equivalent of finitely additive measures for the regular measure $\mu$.


  
{\it Proof.}  
Let us first prove that $(i)\Rightarrow (ii)$.
Let condition $(i)$ be satisfied.
Assume that condition $(ii)$ is false in this case, i.e., there exists $\zeta \in \Delta_{ba}$ such that $\bar{\zeta} \neq \mu$. 
Let us take this measure as the initial measure $\eta=\zeta$ for the Cesaro means from Condition $(i)$.
Then $\lambda_n^\zeta = \zeta$ for all $n\in N$ and for any $f\in C(X)$  will hold
$$ \int_{X}{f(x)\lambda_{n}^{\eta}(dx)} = \int_{X}{f(x)\eta(dx)}= \int_{X}{f(x)\mu(dx)},$$
i.e., $\bar{\eta}=\mu$. We got a contradiction. Therefore, when $(i)$ is executed, $(ii)$ is also executed, i.e., $(i)\Rightarrow (ii)$.

Let us now prove that $(ii)\Rightarrow (i)$.


Let Condition $(ii)$ hold for some $\mu \in S_{rba}$. Assume that $\lambda_{n}^{\eta} \not \to \mu$, in the $\tau_{C}$-topology for some $\eta \in S_{ba}$. Then, by Alexandrov Theorem (see ~\cite{DS1958}, Chapter IY, Section 9, Theorem 15), there exists a open set $G$=$\displaystyle{\buildrel\circ\over G}$ such that $\mu(G)=\mu(\bar{G})$ and $\lambda_{n}^{\eta}(G) \not \to \mu(G)$, i.e., there exist $\varepsilon >0$ and a strictly increasing sequence $\{n_{i}\}$ such that $\lambda_{n_{i}}^{\eta}(G) \ge \mu(G)+\varepsilon$ (or $\le \mu(G)-\varepsilon$) for $i=1,2,\ldots$.

Let $\zeta$ be a $\tau_{B}$-limit point of $\lambda_{n_{i}}^{\eta}$ (here $\tau_{B}$ is the second weak topology we use after $\tau_{С}$). It exists by Theorem 7.2 \cite{Zhd2003a} and Corollary 7.2. \cite{Zhd2003a}. Then $\zeta(G) \ge \mu(G)+\varepsilon $ and, moreover, $\zeta \in \Delta_{ba}$ by Theorem 7.2 \cite{Zhd2003a}.

Since $\bar \zeta$ is regular, for every $\delta > 0$, we can find a set $F = \bar{F} \subset G$ such that $\zeta(F) \ge \zeta(G)-\delta $. The difference $X \setminus G$ is closed and $(X\setminus G)\cap{F} = \varnothing$; therefore, by the Urysohn theorem (see ~\cite{DS1958}, Chapter I, Section 5, Theorem 2), there exists a function $f \in C(X), 0 \le f(x) \le 1$, with $f(F) = 1$ and $f(X\setminus G) = 0$.

Estimate the following integrals:

          $$ \int_{X}{f d\zeta}\ge \int_{F}{f d\zeta}\ge \zeta (F) \ge \zeta (G)-\delta \ge \mu(G)+\varepsilon - \delta $$ 
                               
          $$\ge \int_{G}{f d\mu} +\varepsilon - \delta =  \int_{X}{f d\mu} +\varepsilon - \delta.$$ 

Put $\delta = \frac{\varepsilon}{2}.$   Then 
$$\int_{X}{f d\zeta}\ge \int_{X}{f d\mu} + \frac{\varepsilon}{2},                                
\text{\;\;\; i.e. }, 
\int_{X}{f d \bar{\zeta}} = \int_{X}{f d\zeta} \ne  \int_{X}{f d\mu} , $$ 
and $\bar \zeta \ne \mu$,which contradicts $(ii)$.

Consider the other possible case, $\lambda_{n_i}^\eta(G)\le\mu(G)-\varepsilon$ for $i=1, 2, \ldots$ . Then a $\tau_B$-limit point $\zeta$ of $\lambda_{n_i}^\eta$ satisfies $\zeta(G)\le\mu(G)-\varepsilon$ and $\zeta\in\Delta$. Since $\mu$ is regular, for every $\delta>0$, there exists a set $F=\bar{F}\subset G$ such that $\mu(F)\ge\mu(G)-\delta$, i.e., $\mu(G)\le\mu(F)+\delta$.

Take again a function $f\in C(X)$, $0\le f(X)\le 1$, such that $f(F)=1$ and $f(X\setminus G)=0$. We have
   \begin{eqnarray*}
	\int\limits_X fd\zeta&=&\int\limits_G fd\zeta\le\zeta(G)\le\mu(G)-\varepsilon\le\mu(F)-\varepsilon+\delta=\\
	&=&\int\limits_F fd\mu-\varepsilon+\delta \le \int\limits_X fd\mu-\varepsilon+\delta.
   \end{eqnarray*}
Put $\delta=\frac{\varepsilon}{2}$ and obtain $\bar \zeta \ne \mu$, which contradicts $(ii)$. Consequently, in both cases, $\lambda_{n}^{\eta}\to\mu$ in the $\tau_C$-topology for every $\eta \in S_{ba}$. The theorem is proved. \\

Let us now turn to the case of a topological phase space $(X, \mathfrak B)$. 

{\bf Definition 3.1.} An MC defined on $(X, \mathfrak B)$ is called {Feller} if $TC(X) \subset C(X)$. The Markov operators corresponding to a Feller MC are also called {\it Feller}.
 
If the MC is Feller then $(ii)$ implies $\mu \in \Delta_{rba}$, i.e., $\mu$ is an invariant measure.

Recall that a Feller MC defined on a Hausdorff compact space always has an invariant countably additive probability measure. For a metric compact space, this was proved by Bebutov \cite{Beb1942} (see also \cite{Zhd2003a}, Theorem 4.5).
 
If under the conditions of Theorem 3.1 $\Delta_{ba}$ does not contain purely finitely additive measures then Doeblin condition holds, i.e., the Markov operator $A$ is quasi-compact \cite{Zhd2003b}. 

Note 3.2. 
In the theory of general Markov chains, one usually constructs ergodic theorems for Cesaro means of a Markov sequence of measures $\{\lambda_n^\eta\}$. This is done for the case when it is possible for the MC to have cycles of measures. If the MC has no cycles, i.e., it is acyclic, then in the limit theorems the Cesaro means can be replaced by the Markov sequence of measures $\{\mu_n \}$ itself, which strengthens the ergodicity property of the MC.

Within the framework of our approach, in which the Markov operators $A$ are extended to the space of finitely additive measures, we must also take into account possible cycles consisting of finitely additive measures. Such cycles for general MCs are studied in our paper~\cite{Zhd2021}. 

Let the measure $\mu$ in the conditions of Theorem 3.1 be countably additive. Then, figuratively, invariant purely finitely additive measures are ``buffer'' near the limit countably additive measure $\mu$ (possibly invariant, and possibly being an ``ejection'' point for the operator A). If there is no ``buffer'', then the measures $\lambda_n^\eta$ converges strongly to $\mu$ \cite{Zhd2003a}, \cite{Zhd2003b}. In the presence of a ``buffer'' the MC converges weakly ``sticking'' in invariant purely finitely additive measures ``stuck'' to the limit measure $\mu$ in the $\tau_C$-topology.

However, we can give a completely different interpretation of the assertion from Theorem 3.1. Suppose that conditions $(ii)$ are satisfied and that the Cesaro means $\lambda_n^\eta$ converges $\tau_C$-weakly to the measure $\mu \in S_{rba}$ for any initial measure $\eta \in S_{ba}$. Now, let $\zeta \in \Delta_{ba}$ and $\bar \zeta = \mu$. Since the measure $\zeta$ is not separable from the measure $\mu$ in the $\tau_C$-topology, we have the right to say that the sequence $\lambda_n^\eta$ converges $\tau_C$-weakly to the measure $\zeta$, which is invariant for the MC.

Moreover, this phrase is true irrespective of whether the measure $\mu$ is invariant or not, i.e., for $ \mu \in \Delta_{ba} $ and for $ \mu \notin \Delta_{ba} $.

This same phrase is also true for any other invariant measure $\xi \in \Delta_{ba}$ (there can be infinitely many such measures).

So, we get the following corollary. 

{\bf Corollary 3.1.} Suppose that the conditions of Theorem 3.1 are satisfied. Then for any initial measure $ \eta \in S_{ba} $ the sequence of measures $ \lambda_n^\eta $ $\tau_C$-weakly converges to the measure $ \mu $ and to every invariant measure of the Markov chain $ \zeta \in \Delta_{ba} $, even when the measure $ \mu $ is not invariant.

The first version of Theorem 3.1 was introduced back in 1981 in our paper \cite{Zhd1981} (Theorem 4). There we consider an arbitrary metric space $ (X, \rho) $ on which all countably additive measures are regular. It is assumed that the limit measure $ \mu $ is countably additive (and regular). The condition $(ii)$ for the weak convergence of Cesaro means for MC is the same as in the present paper.

Later (in 2003) in our paper \cite{Zhd2003b} (Theorem 13.1) a generalization of this theorem was proved to the case of an arbitrary normal topological space $ (X, \tau) $, on which countably additive measures may not be regular. Moreover, the limit measure $ \mu $ was assumed to be countably additive and regular.

In this paper, Theorem 3.1, which generalizes Theorem 4 of \cite{Zhd1981} and Theorem 13.1 of \cite{Zhd2003b}, makes an essential weakening of the conditions: the limit measure $ \mu $ can now be only finitely additive (but regular). The scheme of the proof of Theorem 3.1 in this article is similar to the scheme of the proof of Theorem 13.1 in paper \cite{Zhd2003b}. This extends the applicability of Theorem 3.1 to specific MC examples.


\section{\bf Example for using Theorem 3.1.}

Let $X=[0,1]$ be a phase space with Borel $\Sigma$-algebra $\mathfrak B$. On the space $(X, \mathfrak B)$  is given a Markov chain (MC) with a transition function  $p(x,E), \; x\in [0,1], E \in \mathfrak B$ satisfying the following conditions:\medskip

$1. \; \text {If } \; x\in (0,1], E \in \mathfrak B \text { then} $
$$ p(x,E)=\frac{\lambda (E \cap (0, x))}{\lambda((0, x))}=\frac{\lambda(E \cap (0, x))}{x}, $$
$\text {where} \;\; \lambda(E) - \text {Lebesgue measure}.$

Thus, $p(x,E)$ has a uniform distribution on the interval $(0,x)$ for any $x \in (0, 1]$. In particular, $p(1,E)=\lambda(E)$ for any $E \in \mathfrak B$. \medskip

$2. \; \text {If } \; x=0 \; \text { then} $
$$ p(0, \{1\})=1, \;\; p(0, [0,1))=0.  $$ 

It is easy to check that the transition function  $p(x,E)$ is the countably additive  measure for the second argument, i.e., $p(x,E)$ is classical. \medskip

The phaze portrait of our Markov chain is shown in the Figure 1.
\begin{figure}[h!]
	\center{\includegraphics[width=0.9\linewidth]{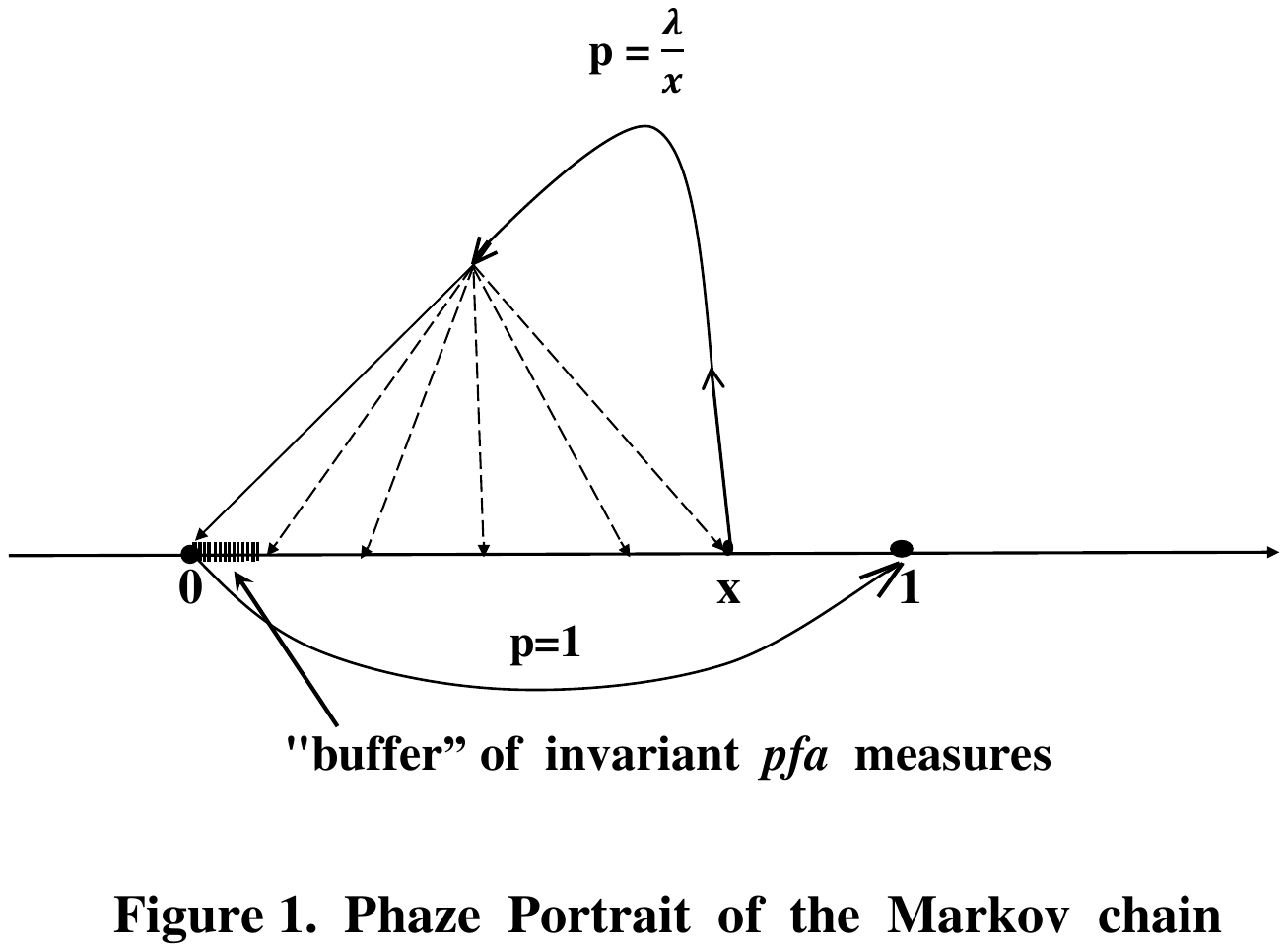}}
\end{figure}  

{Corollary for p(x, E):}

$1) \;\;  p(x, \{z\})=0, \;  \forall x \in (0,1] , \;\; \forall z \in [0,1] ;$

$2) \;\; p(x, (0,x))=1,  \; \forall x \in (0,1);$

$3) \;\; p(x, (0,\varepsilon))=\frac{\varepsilon}{x},  \; \forall \; 0<\varepsilon \le x.$ 

A transition function defines a homogeneous Markov chain (MC) on $([0,1],\mathfrak B)$ with Markov operators $A$ and $T$ (see above).  \medskip


Any Markov chain has an invariant finitely additive measure (\cite{Si1962}) $ \mu = A\mu \in S_{ba}$. We consider its properties for our MC:
$$1) \;\; \mu(\{0\})=A \mu(\{0\})=\int_{[0,1]}p(x,\{0\}) \mu(dx)=\int_{\{0\}} + \int_{(0,1]} =$$
$$= \underbrace{p(0,\{0\})}_{=0}\mu(\{0\}) + \int_{(0,1]}\underbrace{p(x,\{0\})}_{=0} \mu(dx)= 0,$$ \medskip
 i.e., \; $ \mu(\{0\})=0$. 
$$ \;\;\;\;  \;\; \;\; 2) \;\; \mu(\{1\})=A \mu(\{1\})=\int_{[0,1]}p(x,\{1\}) \mu(dx)=\int_{\{0\}} + \int_{\{1\}} + \int_{(0,1)} =$$ 
$$= \underbrace{p(0,\{1\})}_{=1}\mu(\{0\}) + \underbrace{p(1,\{1\})}_{=0}\mu(\{1\}) + \int_{(0,1)}0\cdot \mu(dx)= \mu(\{0\}),$$ 
 i.e., $ \mu(\{1\})=\mu(\{0\})=0$.
 
 $$3) \;\; \text{Let} \;\; 0< \varepsilon < 1.\;\; \text {Then, we have} \;\; \mu((0, \varepsilon))=A \mu((0, \varepsilon))=$$
 $$=\int_{[0,1]}p(x,(0, \varepsilon)) \mu(dx)=\int_{[0, \varepsilon)} + \int_{[\varepsilon,1]} =$$
$$= \int_{[0, \varepsilon)}\underbrace{p(x,(0, \varepsilon))}_{=1} \mu(dx) + \int_{[\varepsilon,1]}p(x,(0, \varepsilon)) \mu(dx) = $$ 
 $$=1 \cdot \mu([0, \varepsilon))+ \int_{[\varepsilon,1]}{\frac{\lambda ((0, \varepsilon))}{\lambda ((0, x))}} \mu(dx)=$$
 $$=\mu([0, \varepsilon))+\int_{[\varepsilon,1]}{\frac{\varepsilon}{x}} \mu(dx)=\mu([0, \varepsilon))+\varepsilon \cdot \int_{[\varepsilon,1]}{\frac{1}{x}} \mu(dx).$$\\
 
 Since $\mu(\{0\})=0$ and $\mu([0,\varepsilon))=\mu((0,\varepsilon))$, we got   
   $$ \mu((0, \varepsilon)) = \mu((0, \varepsilon))+\varepsilon \cdot \int_{[\varepsilon,1]}{\frac{1}{x}} \mu(dx), \;\; \text
{i.e.,} \;\; \varepsilon \cdot \int_{[\varepsilon,1]}{\frac{1}{x}} \mu(dx)=0. $$ \medskip
 
Since $0< \varepsilon \leq x \leq 1$, i.e., $x \in [\varepsilon,1]$, we have ${\frac{1}{\varepsilon}} \geq {\frac{1}{x}} \geq 1$. From here we get 
$$ 0= \int_{[\varepsilon,1]}{\frac{1}{x}} \mu(dx) \geq \int_{[\varepsilon,1]} 1\cdot \mu(dx) = \mu([\varepsilon,1]) \geq 0.$$

Therefore, $\mu([\varepsilon,1]) = 0$ for all $\varepsilon \in (0,1]$.

Whence, $ 1 = \mu([0,\varepsilon)) + \mu([\varepsilon,1]) = \mu([0,\varepsilon))  $ for all $\varepsilon \in (0,1)$.

Since $\mu(\{0\})=0$ then $\mu((0,\varepsilon)) \equiv 1$  for all $\varepsilon \in (0,1)$.


Such measures $\mu$ are typical purely finitely additive measures. \medskip

{\underline {Conclusion}}: our MC does not have invariant countably additive measures but has invariant purely finitely additive measures (such a measure is not unique, see~\cite{Zhd2003a}, Theorem 8.3). All such measures satisfy condition: $\mu((0,\varepsilon)) = 1$ and $\mu([\varepsilon, 1)) = 0$  for all $\varepsilon \in (0,1)$,  $\mu(\{0\})=0$. So, $\Delta_{ca} = \varnothing, \Delta_{pfa} \ne \varnothing$. \medskip

We verify the fulfillment of condition $(ii)$ of Theorem 3.1.

 If $f \in C[0,1]$ then for any invariant measure $\mu \in S_{ba}$ of our MC and for arbitrary $\varepsilon > 0$ and we have 
$$ \int_{[0,1]}{f(x)\mu(dx)} = \int_{[0,\varepsilon]}{f(x)\mu(dx)} = \lim\limits_{x \to 0} f(x) = f(0).$$

Let $\delta_{0}$ be the Dirac measure degenerate at the point $0$, i.e., $\delta_{0}(\{0\})=1, \;\; \delta_{0}((0,1])=0$.
The measure $\delta_{0}$ is a regular countably additive measure.

Then for any $f \in C[0,1]$,  $\int_{[0,1]}{f(x)\delta_{0}(dx)} = f(0)$.

Therefore, the condition $(ii)$ of Theorem 3.1 is satisfied: $$\int_{[0,1]}{f(x) d\mu(x)} = \int_{[0,1]}{f(x)\delta_{0}(dx)}$$ for any $f \in C[0,1]$ and for any invariant finitely additive measure $\mu$ of our MC.

Now, as was proved in Theorem 3.1, i.e., condition $(i)$ is satisfied, for any initial measure $\eta \in S_{ba}$ the Cesaro means $\lambda_{n}^{\eta}$ for a Markov chain $\tau_{C}$-weakly converge to the limiting measure $\delta_{0}$, i.e., for all $ f \in C[0,1]$ we have   
$$ \int_{[0,1]}{f(x)\lambda_{n}^{\eta}(dx)} \to \int_{[0,1]}{f(x)\delta_{0}(dx)}.$$

{Here a countably additive regular limit measure $\bf \delta_{0}$ is not invariant for MC, i.e., $ {A\delta_{0} \ne \delta_{0}}$.} Really, $A\delta_{0}(\{0\}) = p(0,\{0\})=0$,  $A\delta_{0}(\{1\}) = p(0,\{1\})=1$, i.e.,  $A\delta_{0} =  \delta_{1}$, where $\delta_{1}$ is the Dirac measure degenerate at the point $\{1\}$. And we have $A\delta_{1}=A^{2}\delta_{0} = \lambda$ in the next step.

It can be shown that the sequence $\lambda_{n}^{\eta}$ for any countably additive measures $\eta$ does not converge to $\delta_{0}$ in any of the topologies $\tau_{ba}$(strong), $\tau_{ba^{*}}$(weak), $\tau_{B}$(*-weak).\medskip

In accordance with Corollary 3.1, we have the right to say that the Cesaro means $\lambda_{n}^{\eta}$ for any initial measure $\eta \in S_{ba}$ simultaneously converge $\tau_{C}$-weakly not only to the measure $\mu=\delta_{0}$, but also to {\it each} invariant purely finitely additive measure $\xi \in \Delta_{pfa}$ that is not separable from $\mu$ in the $\tau_{C}$-topology and such a measure $\xi$ is not unique.\medskip

For completeness of the description of this MC, we also consider the traditional question of its Feller property, since within the framework of our operator approach, the Feller property of MC is related to invariant measures.

Our MC is given on the metric compactum $ X = [0,1] $. If it were Feller, then, according to the well-known theorem of Bebutov~\cite{Beb1942}, it would have an invariant bounded countably additive measure. And there are not such invariant measures, as we showed above. Therefore, this MC is not a Feller chain.

However, it is interesting to know exactly where $ X = [0, 1] $ is already violated by Feller's (we intuitively assume that at $ x = 0 $). Let us find these points.\medskip

Let $f \in C_{[0,1]}, \; g(x) = Tf(x) = \int_{[0,1]} f(y) p(x,dy), \; x \in [0,1]$.
	
Let us consider some cases.

1) $g(0) = Tf(0) = \int_{[0,1]} f(y) p(0,dy) = f(1) p(0,\{1\}) = f(1)$.

2) $g(1) = Tf(1) = \int_{[0,1]} f(y) p(1,dy) = \int_{[0,1]} f(y) \lambda(dy) = (R) \int_{[0,1]} f(y)dy$. (Here and below (R) denotes the Riemann integral).

3) Let $0<x<1$. Then $$g(x) = Tf(x) = \int_{[0,1]} f(y) p(x,dy) = \int_{[0,x)}f(y)\cdot\frac{1}{x}\cdot\lambda(dy) =$$
$$ \frac{1}{x} \cdot \int_{[0,x)} f(y)  \lambda(dy) = \frac{1}{x}(R)\int_{0}^{x} f(y)dy.$$ 

Since the definite Riemann integral of a continuous function is continuous with respect to the upper limit, and the function $ \frac {1}{x} $ is continuous for $ x \ne 0 $, the function $ g(x) $ is also continuous for any $ x \in (0,1) $.\medskip

We now verify the continuity of the function $ g = Tf $ for $ x = 1 $ and $ x = 0 $.

4) $$\lim\limits_{x \to 1} g(x) = \lim\limits_{x \to 1} (\frac{1}{x} \cdot \int_{0}^{x} f(y)dy) = \int_{0}^{1} f(y)dy = g(1).$$ 

Consequently, the function $ g (x) $ is continuous at the point $x=1$. 

5)  $$\lim\limits_{x \to 0} g(x) = \lim\limits_{x \to 0} (\frac{1}{x} \cdot \int_{0}^{x} f(y)dy).$$

Let us show that, generally speaking, $\lim\limits_{x \to 0} g(x) \ne g(0) = f(1)$.

We put $f(y) = y \in [0,1], f \in C_{[0,1]}$. Then $$g(x) = \frac{1}{x} \cdot \int_{0}^{x} f(y)dy = \frac{1}{x} \cdot \int_{0}^{x} ydy = \frac{x^2}{2\cdot x} = \frac{x}{2},$$ and $\lim\limits_{x \to 0} g(x) = \lim\limits_{x \to 0} \frac{x}{2} =0.$ 

But $g(0) = f(1) =1$, i.e., $\lim\limits_{x \to 0} g(x) \ne g(0)$.

Consequently, the function $ g(x) $ has a discontinuity at the point $x=0$, and $g \notin C_{[0,1]}$. 

So, we showed that our MC is really not Feller on the compact $ [0,1] $. We call such MC are almost Feller.

Now it becomes clear why the buffer from all invariant purely finitely additive measures of our MC accumulated around the point $ 0 $.

Note 4.1.
All the review of this example was carried out without using the explicit form of measures $ \mu_{n}=A^{n}\eta, \lambda_{n}^{\eta} $ and  $\int{f(x)\lambda_{n}^{\eta}(dx)}.$
An explicit description of these measures for various $\eta$ is a not simple analytical question.

Our Theorem 3.1 made it possible to obtain information on the limiting behavior of MC by using only qualitative statements on invariant purely finitely additive measures for MC, which turned out to be much simpler.\medskip

Note 4.2. 
In the theory of dynamical systems, ergodic theorems are also studied in a number of papers in the absence of invariant countably additive measures. However, these are not exactly measures (and in some cases not at all) that are called invariant in the theory of Markov chains. But, there are points (and areas) of contact between these two concepts and theories (see, for example, our work \cite{Zhd1982}). In this paper we do not consider these problems.



\begin{thebibliography}{99}
	
\bibitem{YoKa:1}  
\emph{Yosida~K., Kakutani~S.} Operator-theoretical treatment of Markoff's processes and mean ergodic theorem.~Ann. Math., 1941, V 42, №1. 188-228. 

\bibitem{Zhd84} 
\emph{Zhdanok A.I.}	Regularization of finitely additive measures.~Latv. Mat. Ezegodnik, v. 28, 1984. 234-248, Zinatne, Riga (Russian).

\bibitem{Zhd2003a} 
\emph{Zhdanok A.I.} Finitely additive measures in the ergodic theory of Markov chains. I.~{\it Siberian Advances in Math.}, 2003, v. 13, N 1. 87--125. 

\bibitem{Zhd2003b} 
\emph{Zhdanok A.I.} Finitely additive measures in the ergodic theory of Markov chains. II.~{\it Siberian Advances in Math.}, 2003, v. 13, N 2. 108--125. 

\bibitem{Zhd1981} 
\emph{Zhdanok A.I.}	Invariant finitely additive measures and the limit behavior of Markov processes with discrete time.~{\it Dokl. Akad. Nauk Ukrain. SSR Ser. A}, 1981, N3. 11-13 (Russian).

\bibitem{DS1958} 
\emph{Dunford ~N, Schwartz ~J.} Linear operatiors. Part I: General Theory. ~-Interscience Publischer, New York, London,~--1958. (Данфорд~Н., Щварц~Дж. Линейные операторы. Общая теория.М., ИЛ, 1962. 896с.) 

\bibitem{YoHew1952} 
\emph{Yosida~K., Hewitt~E.} Finitely additive measures.~Trans. Amer. Math. Soc., 1952, v. 72, I. 46--66. 

\bibitem{Rama1981} 
\emph{Ramakrishnan S.} Finitely Additive Markov Chains.~Trans. Amer. Math. Soc., 1981, v. 265. 247--272. 


\bibitem{ZhdKh2022} 
\emph{Zhdanok A., Khuruma A.} Decomposition of Finitely Additive Markov Chains in Discrete Space.~ Mathematics, v. 10(12),2083, 2022. 21 p. 

\bibitem{Si1962} 
\emph{\v{S}idak~Z}. Integral representations for transition probabilities of Markov chains with a general state space.~Czechoslovak Math. J., 12 (87), 4,  1962. 492--522. 

\bibitem{Fog1966} 
\emph{Foguel S.~R.} Existence of Invariant Measures for Markov Processes. II.~
Proc. Amer. Math. Soc., 17, no. 2, 1966. 387--389.

\bibitem{Beb1942} 
\emph{Bebutov M.B.} Markov chains with compact state space.~Mat. Sb., N10, 1942. 231--238 (Russian).

\bibitem{Zhd2021} 
\emph{Zhdanok A.I.}	Cycles in Spaces of Finitely Additive Measures of General Markov Chains.~Recent Developments in Stochastic Methods and Applications. Springer Proceedings in Mathematics  Statistics, v. 371, 2021. 153--164.

\bibitem{Zhd1982} 
\emph{Zhdanok A.I.}	Iterative processes as Markov chains.~Latv. Mat. Ezegodnik, v. 26, 1982. 153--164, Zinatne, Riga (Russian).

\end{thebibliography}
\end{document}